\documentclass[a4paper,11pt]{amsart}
\usepackage{amssymb}
\usepackage{ifthen}
\usepackage{graphicx}
\nonstopmode
\newtheorem{thm}{Theorem}
\newtheorem{cor}{Corollary}
\newtheorem{lem}{Lemma}

\newtheorem{claim}{Claim}

\newtheorem{conj}{Conjecture}
\newtheorem{prob}{Problem}

\theoremstyle{definition}
\newtheorem{defn}{Definition}
\newtheorem{example}{Example}

\newenvironment{rem}{%
\bigskip
\noindent \textsl{{\sl Remark. }}}{\bigskip}
\newenvironment{rems}{%
\bigskip
\noindent \textsl{{\sl Remarks. }}}{\bigskip}

\newenvironment{pf}[1][]{%
 \vskip 1mm
 \noindent
 \ifthenelse{\equal{#1}{}}%
  {{\slshape Proof. }}%
  {{\slshape #1.} }%
 }%
{\qed\bigskip}

\newcounter{alphabet}
\newcounter{tmp}

\makeatletter
\newcommand{\Ref}[1]{\@ifundefined{r@#1}{}{\setcounter{tmp}{\ref{#1}}\Alph{tmp}}}
\makeatother

\newenvironment{Lem}[1][]{\refstepcounter{alphabet}%
\bigskip%
\noindent%
{\bf Lemma \Alph{alphabet}}%
{\bf .} \itshape}{\vskip 8pt}

\newcommand{\IC}{{\mathbb C}}
\newcommand{\ID}{{\mathbb D}}

\def\be{\begin{equation}}
\def\ee{\end{equation}}

\newcommand{\bee}{\begin{enumerate}}
\newcommand{\eee}{\end{enumerate}}

\newcommand{\blem}{\begin{lem}}
\newcommand{\elem}{\end{lem}}
\newcommand{\bthm}{\begin{thm}}
\newcommand{\ethm}{\end{thm}}
\newcommand{\bcor}{\begin{cor}}
\newcommand{\ecor}{\end{cor}}
\newcommand{\beg}{\begin{example}}
\newcommand{\eeg}{\end{example}}
\newcommand{\begs}{\begin{examples}}
\newcommand{\eegs}{\end{examples}}
\newcommand{\bdefe}{\begin{defn}}
\newcommand{\edefe}{\end{defn}}
\newcommand{\bprob}{\begin{prob}}
\newcommand{\eprob}{\end{prob}}
\newcommand{\bques}{\begin{ques}}
\newcommand{\eques}{\end{ques}}
\newcommand{\bei}{\begin{itemize}}
\newcommand{\eei}{\end{itemize}}

\newcommand{\bde}{\begin{deter}}
\newcommand{\ede}{\end{deter}}
\newcommand{\bca}{\begin{case}}
\newcommand{\eca}{\end{case}}
\newcommand{\bcl}{\begin{claim}}
\newcommand{\ecl}{\end{claim}}
\newcommand{\bcon}{\begin{conj}}
\newcommand{\econ}{\end{conj}}
\newcommand{\bcons}{\begin{conjs}}
\newcommand{\econs}{\end{conjs}}
\newcommand{\bprop}{\begin{propo}}
\newcommand{\eprop}{\end{propo}}
\newcommand{\br}{\begin{rem}}
\newcommand{\er}{\end{rem}}
\newcommand{\brs}{\begin{rems}}
\newcommand{\ers}{\end{rems}}
\newcommand{\bo}{\begin{obser}}
\newcommand{\eo}{\end{obser}}
\newcommand{\bos}{\begin{obsers}}
\newcommand{\eos}{\end{obsers}}
\newcommand{\bpf}{\begin{pf}}
\newcommand{\epf}{\end{pf}}
\newcommand{\ba}{\begin{array}}
\newcommand{\ea}{\end{array}}
\newcommand{\beq}{\begin{eqnarray}}
\newcommand{\beqq}{\begin{eqnarray*}}
\newcommand{\eeq}{\end{eqnarray}}
\newcommand{\eeqq}{\end{eqnarray*}}

\newcommand{\ds}{\displaystyle}

\begin{document}
\title[On the generalized Zalcman functional $\lambda a_n^2-a_{2n-1}$ in the close-to-convex
family]{On the generalized Zalcman functional $\lambda a_n^2-a_{2n-1}$ in the close-to-convex family}


\author[L. Li]{Liulan Li}
\address{Liulan Li, College of Mathematics and Statistics,
Hengyang Normal University, Hengyang,  Hunan 421002, People's
Republic of China.} \email{lanlimail2012@sina.cn}

\author[S. Ponnusamy]{Saminathan Ponnusamy
}
\address{Saminathan  Ponnusamy,
Indian Statistical Institute (ISI), Chennai Centre, SETS (Society
for Electronic Transactions and Security), MGR Knowledge City, CIT
Campus, Taramani, Chennai 600 113, India.}
\email{samy@isichennai.res.in, samy@iitm.ac.in}



\begin{abstract}
Let  ${\mathcal S}$ denote the class of all functions $f(z)=z+\sum_{n=2}^{\infty}a_{n}z^{n}$ analytic and univalent in the unit disk $\ID$.
For $f\in {\mathcal S}$, Zalcman conjectured that
$|a_n^2-a_{2n-1}|\leq (n-1)^2$ for $n\geq 3$. This conjecture has been verified only certain values of $n$ for $f\in {\mathcal S}$ and for all
$n\ge 4$ for the class $\mathcal C$ of close-to-convex functions (and also for a couple of other classes). In this paper we provide bounds of the
generalized Zalcman coefficient functional $|\lambda a_n^2-a_{2n-1}|$ for functions in $\mathcal C$ and for all $n\ge 3$,
where $\lambda$ is a positive constant. In particular, our special case settles the open problem
on the Zalcman inequality for $f\in \mathcal C$ (i.e. for the case $\lambda =1$ and $n=3$).
\end{abstract}

\maketitle \pagestyle{myheadings} \markboth{L. Li and S.
Ponnusamy}{Generalized Zalcman functional in the close-to-convex family}

\section{Introduction and Main results}

Let $\mathbb{D}=\{z:\,|z|<1\}$ be the open unit disk in $\IC$ and $\partial \mathbb{D}=\{z:\,|z|=1\}$. Let
${\mathcal A}$ denote the class of all analytic functions in $\mathbb{D}$ and ${\mathcal S}\subset {\mathcal A}$
denote the family of all normalized univalent functions $f$ of the form
\be\label{eq1}
f(z)=z+\sum_{n=2}^{\infty}a_{n}z^{n}.
\ee
A function $f\in {\mathcal S}$ is called starlike if $f(\ID)$ is starlike with respect to the origin.
Let  ${\mathcal S}^*$ denote the class of all starlike functions in $\mathcal S$.
An analytic function $f$ of the form \eqref{eq1} is called close-to-convex if there exists a real number $\theta$ and a function $g\in {\mathcal S}^*$
such that $ {\rm Re\,} \left(e^{i\theta}zf'(z)/g(z)\right)>0$ for $z\in \ID$. Functions in the class  ${\mathcal C}$  of all close-to-convex functions
are known to be univalent in $\ID$. Geometrically, $f\in {\mathcal C}$ means that the complement of the
image-domain $f(\ID)$ is the union of rays that are disjoint (except that the origin of one ray may lie on another one of the rays).

The role of the family ${\mathcal S}$ together with their subfamilies and their importance  geometric function theory
are well documented, see for example \cite{Du,Go,HM,Pom,Pom2}. Several conjectures which imply the Bieberbach conjecture that
$|a_n|\leq n$ for each $f\in{\mathcal S}$ have been verified by the de Branges theorem.
At that time, as an approach to prove the Bieberbach conjecture (see also \cite{BT}), Lawrence Zalcman in 1960 conjectured that the coefficients
of ${\mathcal S}$ satisfy the sharp inequality
\be\label{eq1a}
|a_n^2-a_{2n-1}|\leq (n-1)^2,
\ee
for each $n\geq 2$ with equality only for the Koebe function $k(z)=z/(1-z)^2$ and its rotation.
The conjecture remains open and partial progress was obtained in the sense that
it has been verified either for certain  subclasses of ${\mathcal S}$ or for certain values of $n$ for the full class ${\mathcal S}$
(eg. see \cite{Kr95} for the case $n=3$,  and  \cite{Kr10} for the cases $n=4,5,6$). Analog of this conjecture for several other subclasses of
${\mathcal S}$ has not been attempted except for some cases \cite{AbLiPo-2014,BhPoWir2011,LiPo-2014,LiPo-2015}.

Recall that the case $n=2$ of Zalcman conjecture is the well-known
Fekete-Szeg\"{o} inequality, namely, $ |a_2^2-a_{3}|\leq 1 $ (see \cite[Theorem 1.5]{Pom} and \cite[Theorem 3.8]{Du}).
Moreover, the problem of maximizing $|\lambda a_2^2-a_{3}|$ with $\lambda >0$ for $\mathcal S$ (\cite[Theorem 3.8]{Du}) and for
many other subclasses has been solved by a number of authors--often for restricted values of  $\lambda$, see for instance
\cite{BhPoWir2011,CKSugawa2010, KanasDar2010,KCSugawa2000,Koepf87,Koepf87b}. In \cite{Pflu-86},  Pfluger extended the later result
for certain complex values of $\lambda$.

In \cite{Pflu-76}, Pfluger pointed out that if  $f\in{\mathcal S}$ then the coefficients of
$\sqrt{f(z^2)}$ and $1/f(1/z)$ are polynomials in $a_j$ which contain expressions of the form $\lambda a_n^2-a_{2n-1}$ which is referred to as the
generalized Zalcman functional. Thus it is natural to consider the problem of maximizing $|\lambda a_n^2-a_{2n-1}|$, as a function of the real parameter $\lambda$.
There are only few results of this type in the literature and that too deal with only $\lambda =1$ for a few  classes of functions $f$. For example,
in \cite{BT}, Brown and Tsao proved the Zalcman conjecture for $n\geq 3$ for typically real functions (hence for functions with real coefficients in $\mathcal S$)
and also for the class ${\mathcal S}^*$.
In fact they have proved it in a general form involving generalized Zalcman functional, and another generalization of the result of
Brown and Tsao appeared in \cite{MaWi}.
In 1988,  Ma \cite{MaWa} provided further evidence in support of the conjecture by verifying the Zalcman inequality \eqref{eq1a} (again for the case $\lambda =1$ only)
for the class ${\mathcal C}$ when $n\geq 4$. However, the conjecture remains open in ${\mathcal C}$ for $n=3$ (see \cite[Remarks]{MaWa}).

In this article, we consider the generalized Zalcman conjecture for the class ${\mathcal C}$ and $n\geq3$.
In particular, we solve the conjecture for the family $\mathcal C$ when $n=3$. Also, one of our main results contains the proof of the main
result of Ma \cite{MaWa}.

We end the section by indicating the recent result on generalized Zalcman conjecture. Let ${\mathcal F}(\alpha )$ denote the family
of convex functions of order $\alpha$ with  $-1/2\leq \alpha <1$, consisting of functions $f\in  {\mathcal S}$  satisfying the condition
$${\rm Re } \left( 1+\frac{zf''(z)}{f'(z)}\right) > \alpha   ~\mbox{ for $z\in \ID$}.
$$
We remark that ${\mathcal F}(0):={\mathcal K}$, the class of all normalized convex univalent functions, i.e.  each $f\in {\mathcal K}$
maps $\ID$ univalently onto a convex domain. It is known that
${\mathcal F}(\alpha )\subset {\mathcal C}$ for  $-1/2\leq \alpha <0$, and ${\mathcal F}(\alpha )\subset {\mathcal K}$ for  $0\leq \alpha <1$.
The sharp bound for the quantity $|a_n^2-a_{2n-1}|$ for the class of convex functions of order $-1/2$ was first discussed in \cite{AbLiPo-2014}
and in recent articles, Li and Ponnusamy \cite{LiPo-2014} (see also Li et al.  \cite{LiPo-2015}) obtained sharp bound for the
generalized Zalcman functional $|\lambda a_n^2-a_{2n-1}|$, for the family ${\mathcal F}(\alpha )$.

We now state our main results. It is worth pointing out that all the results below continue to hold if ``$f\in {\mathcal C}$" is replaced by
``$f\in {\mathcal S}^*$". This observation shows that our results include a proof of the generalized Zalcman conjecture for the class ${\mathcal S}^*$
and $n\geq3$.

\begin{thm}\label{main theorem}
Suppose that $f\in {\mathcal C}$ as in the form \eqref{eq1} and $n\geq 3$.
\bee
\item [{\rm (1)}] If $\lambda\geq\frac{3n+\sqrt{5n^2-4n}}{n^2+n},$ then
$$|\lambda a_n^2-a_{2n-1}|\leq \lambda n^2-(2n-1),
$$
where the equality holds if $f(z)$ is  the Koebe function $z/(1-z)^2$.

\item [{\rm (2)}] If $\frac{2n}{n^2-n+1}<\lambda<\frac{3n+\sqrt{5n^2-4n}}{n^2+n}$, then
$$|\lambda a_n^2-a_{2n-1}|\leq \frac{\lambda\left[4n(n+1)-(3n^2+1)\lambda\right]+4n(2-\lambda)\sqrt{n\lambda(2-\lambda)}}
{\lambda\left[8n-\lambda(n+1)^2\right]}.
$$

\item [{\rm (3)}] If $0<\lambda\leq\frac{2n}{n^2-n+1}$, then we have
$$|\lambda a_n^2-a_{2n-1}|\leq 2n-1.
$$
\eee
\end{thm}

Substituting $\lambda =1$ in Theorem \ref{main theorem}(1), it follows easily that
$|a_n^2-a_{2n-1}|\leq (n-1)^2$
for all $n\geq 4$ and  $f\in {\mathcal C}$. It is worth to state the cases $n=3, 4$ explicitly.

\begin{cor}\label{n=4}
Suppose that $f\in {\mathcal C}$ as in the form \eqref{eq1}.
\bee
\item [{\rm (1)}] If $\lambda\geq1,$ then
$$|\lambda a_4^2-a_7|\leq 16\lambda -7,
$$
where the equality holds if $f(z)$ is  the Koebe function $z/(1-z)^2$.

\item [{\rm (2)}] If $\frac{13}{8}<\lambda<1$, then
$$|\lambda a_4^2-a_7|\leq \frac{\lambda(80-49\lambda)+32(2-\lambda)\sqrt{\lambda(2-\lambda)}} {\lambda(32-25\lambda)}.
$$

\item [{\rm (3)}] If $0<\lambda\leq\frac{8}{13}$, then we have $|\lambda a_4^2-a_7|\leq 7.$
\eee
\end{cor}

\begin{cor}\label{n=3}
Suppose that $f\in {\mathcal C}$ as in the form \eqref{eq1}.
\bee
\item [{\rm (1)}] If $\lambda\geq\frac{9+\sqrt{33}}{12},$ then
$$|\lambda a_3^2-a_5|\leq 9\lambda -5,
$$
where the equality holds if $f(z)$ is  the Koebe function $z/(1-z)^2$.

\item [{\rm (2)}] If $\frac{6}{7}<\lambda<\frac{9+\sqrt{33}}{12}$, then
$$|\lambda a_3^2-a_5|\leq \frac{\lambda(12-7\lambda)+3(2-\lambda)\sqrt{3\lambda(2-\lambda)}}
{\lambda(6-4\lambda)}.
$$

\item [{\rm (3)}] If $0<\lambda\leq\frac{6}{7}$, then we have $|\lambda a_3^2-a_5|\leq 5.$
\eee
\end{cor}

Clearly, Corollary \ref{n=3}(2) gives that if $f\in {\mathcal C}$ is given by \eqref{eq1}, then
$$| a_3^2-a_5|\leq \frac{5+3\sqrt{3}} {2} \approx 5. 098.
$$


Proof of Theorem \ref{main theorem} rely on a number of lemmas. In Section \ref{sec2}, we present
three important lemmas which play vital role in the formulation of several lemmas in Section \ref{sec3}.   In Section \ref{sec3},
we state and prove several lemmas based on different interval range values of $\lambda$.
The proof of Theorem \ref{main theorem} will be given in Section \ref{sec4}.

\section{Preliminaries and some basic lemmas}\label{sec2}

Suppose that $X$ is a linear topological space and that $Y\subset X$. The subset
$Y$ is called convex if $tx+(1-t)y\in Y$ whenever $x,\ y\in Y$
and $0\leq t\leq 1$. The closed convex hull of $Y$ is defined as the
intersection of all closed convex sets containing $Y$. A point $u\in
Y$ is called an extremal point of $Y$ if, for $0<t<1$ and $x,\ y\in Y$, $u=tx+(1-t)y$
implies that $x=y$. The set $\mathcal{E}Y$
consists of all the extremal points of $Y$ (see \cite{HM, MacGreW}
for a general reference and for many important results on this
topic).

\begin{Lem}\label{extreme} \cite{HM}
Let $\mathcal{HC}$ and $\mathcal{EHC}$ denote the closed convex hull of $\mathcal{C}$ and the set of the extremal points of
$\mathcal{HC}$, respectively. Then $\mathcal{HC}$ consists of all analytic functions represented by
$$ f(z)=\int_{S} \frac{z-\frac{1}{2}(x+y)z^2}{(1-y z)^2}d\mu(x,y),
$$
where $\mu$ is a probability measure on $S=\partial\mathbb{D}\times\partial\mathbb{D}$. The set $\mathcal{EHC}$ consists of
the functions given by
\be\label{eq2}
f(z)=\frac{z-\frac{1}{2}(x+y)z^2}{(1-y z)^2}=z+\sum^\infty_{n=2}\left(\frac{n+1}{2}y^{n-1}-\frac{n-1}{2}x y^{n-2}\right)z^n,
\ee
where $|x|=|y|=1$ and $x\neq y$.
\end{Lem}
\medskip

If the family $\mathcal{F}\subset\mathcal{A}$ is convex and $L:\,\mathcal{A}\rightarrow \mathbb{R}$ is a real-valued functional on
$\mathcal{A}$, then we say that $L$ is convex on $\mathcal{F}$ provided that
$$L\left(tg_1+(1-t)g_2\right)\leq t L(g_1)+(1-t)L(g_2)
$$
whenever $g_1,\ g_2\in\mathcal{F}$ and $0\leq t\leq1.$ Since $\mathcal{HC}$ is convex, we have a real-valued,
continuous and convex functional on $\mathcal{HC}$.

\begin{lem}\label{convex} For $g(z)=z+\sum_{n=2}^{\infty}b_{n}z^{n}$ is analytic in $\mathbb{D}$, consider
$$J(g)=\lambda \left({\rm Re\,} b_{n}\right)^2-{\rm Re\, } b_{2n-1},
$$
where $\lambda>0$. Then $J$ is a real-valued, continuous and convex functional on
$\mathcal{HC}$.
\end{lem}
\bpf Let $h(z)=z+\sum_{n=2}^{\infty}c_{n}z^{n}$ be analytic in
$\mathbb{D}$,  $0\leq t\leq 1$ and $G=tg+(1-t)h$.
Then, by the definition of $J$, we have
$$J(G)= \lambda \left[ {\rm Re\,} (tb_{n} +(1-t)c_n)\right]^2- {\rm Re\,} [tb_{2n-1} +(1-t)c_{2n-1}]
$$
which may be easily rearranged as
$$J(G)= t J(g)+(1-t)J(h) -\lambda t (1-t)[{\rm Re\,} b_{n} - {\rm Re\,} c_n ]^2.
$$
This gives $J(tg+(1-t)h) \leq  t J(g)+(1-t)J(h)$ and the desired conclusions follow.
\epf

Since $\mathcal{C}$ is compact, for the functional $J$ defined as in Lemma \ref{convex}, Theorem~4.6 in \cite{HM} yields
the following.

\begin{Lem}\label{maximum}
$\ds\max\{J(f):\, f\in\mathcal{HC}\}=\max\{J(f):\, f\in\mathcal{C}
\}=\max\{J(f):\, f\in\mathcal{EHC}\}.$
\end{Lem}

By using Lemmas \Ref{extreme}, \ref{convex} and \Ref{maximum}, we derive the following lemma.

\begin{lem}\label{bound}
We have $4\max\{J(f):\ f\in\mathcal{C}\}-4n\leq F_{n, \lambda}(u,v),$
where $F_{n,\lambda}(u,v)=:F(u,v)$ is given by
$$F(u,v)=[(n+1)^2\lambda-8n]u^2-2(n-1)\left[(n+1)\lambda-2\right]uv+4(n-1)\sqrt{1-u^2}\sqrt{1-v^2}+(n-1)^2\lambda v^2,
$$
and $(u,v)\in R=:\,[-1,1]\times [-1,1]$.
\end{lem}
\bpf Let $f(z)=z+\sum_{n=2}^{\infty}a_{n}z^{n}\in\mathcal{EHC}$. Then Lemma \Ref{extreme} and \eqref{eq2}
(with $x=e^{is}$ and  $y=e^{it}$) give that
$$a_n=\frac{n+1}{2}e^{i(n-1)t}-\frac{n-1}{2}e^{i[s+(n-2)t]},
$$
where $t, s\in[0, 2\pi)$ and $t\neq s$. This representation quickly yields
\beqq
J(f)&=&\lambda\left\{\frac{n+1}{2}\cos(n-1)t-\frac{n-1}{2}\cos (s+(n-2)t)\right\}^2 \\
& & -n\cos2(n-1)t+(n-1)\cos (s+(2n-3)t ).
\eeqq
By the identity $\cos 2\theta =2\cos^2\theta -1$ and the addition formula for $\cos (s+(n-2)t +(n-1)t)$, we may rewrite the last expression
as
\beqq
J(f) &=&\left (\frac{(n+1)^2}{4}\lambda-2n\right)\cos^2(n-1)t\\
& & -(n-1)\left (\frac{n+1}{2}\lambda-1\right)\cos(n-1)t\cos (s+(n-2)t)\\
& &~~~ -(n-1)\sin(n-1)t\sin (s+(n-2)t )+\frac{(n-1)^2}{4}\lambda\cos^2 (s+(n-2)t)+n.
\eeqq
If we set $\cos(n-1)t=u$ and $\cos (s+(n-2)t)=v$ and use the identity $\sin\theta=\pm\sqrt{1-\cos^2\theta}$, then the above equation reduces to
$$
4J(f)-4n\leq F(u,v),
$$
where $u,v\in[-1,1]$. Lemmas \Ref{extreme}, \ref{convex} and \Ref{maximum} show that the proof is completed.
\epf

\section{Main lemmas}\label{sec3}

Throughout  $F(u,v):=F_{n, \lambda}(u,v)$ is given by Lemma \ref{bound} whereas $G_{n,\lambda}(u,v)=:G(u,v)$ is defined by \eqref{eq5} below,
and  $R=[-1,1]\times [-1,1]$. Observe that $F(u,v)=F(-u,-v)$ and $G(u,v)=G(-u,-v)$.



\begin{lem}\label{real function bound}
Suppose that  $n\geq 3$ and $\lambda\geq\frac{10n-2}{(n+1)^2}$. Then
$$F(u,v)\leq 4\lambda n^2-12n+4 ~\mbox{ for }~ (u,v)\in R,
$$
where the equality holds if and only if $(u,v)=(1,-1)$ or $(u,v)=(-1,1)$.
\end{lem}
\bpf The elementary inequality $2ab\leq a^2+b^2$ gives
\be\label{eq4}
F(u,v)\leq G(u,v) +4(n-1),
\ee
where $G_{n,\lambda}(u,v)=:G(u,v)$ is given by
\be\label{eq5}
G(u,v)=[\lambda(n+1)^2-10n+2]u^2-2(n-1)[\lambda(n+1)-2]uv+(n-1)[\lambda(n-1)-2]v^2.
\ee
We now introduce
$$A(n)=\frac{10n-2}{(n+1)^2} ~\mbox{ and }~  B(n)=\frac{2}{n}.
$$
Then for $n\geq 3$, we observe that $A(n)>B(n-1)>B(n+1),$
and the functions $A(n)$ and $B(n)$ are monotonically decreasing with respect to $n$. Using $A(n)$ and $B(n)$,
we may write \eqref{eq5} in a convenient form as
\be\label{eq6}
G(u,v)=(n+1)^2\left[\lambda-A(n)\right]u^2-2(n^2-1)\left[\lambda-B(n+1)\right]uv+(n-1)^2\left[\lambda-B(n-1)\right]v^2.
\ee
If $n\geq 3$ and $\lambda\geq\frac{10n-2}{(n+1)^2}$, then $\lambda\geq A(n)>B(n-1)>B(n+1)$ and thus,
$$G(u,v)\leq 4\lambda n^2-16n+8,
$$
where the equality holds if and only if $(u,v)=(1,-1)$ or $(u,v)=(-1,1)$. The proof is completed.
\epf

\begin{lem}\label{G boundary bound 1}
Suppose that  $n\geq 3$ and $\frac{6n-2}{n^2+n}\leq\lambda<\frac{10n-2}{(n+1)^2}$. Then
$$F(u,v)\leq 4\lambda n^2-12n+4~\mbox{ for }~ (u,v)\in R,
$$
where the equality holds if and only if $(u,v)=(1,-1)$ or $(u,v)=(-1,1)$.
\end{lem}
\bpf In view of the inequality \eqref{eq4}, it suffices to prove that
$$G(u,v)\leq G(1,-1)=G(-1,1)=4\lambda n^2-16n+8.
$$

To complete the proof, by assumption, we begin to observe that
\be\label{eq7a}
\frac{2}{n+1}<\frac{2n}{n^2-n+1}<\frac{4n-2}{n^2}<\frac{6n-2}{n^2+n}\leq\lambda<\frac{10n-2}{(n+1)^2}.
\ee
Next, we determine the critical points of $G(u,v)$ and for that we need to consider the equations
\be\label{eq-ex1}
\left \{ \begin{array}{rl}
\ds \frac{1}{2}\cdot\frac{\partial G(u,v)}{\partial u}
=\left[\lambda(n+1)^2-10n+2\right]u-(n-1)\left[\lambda(n+1)-2\right]v=0\\
 \ds \frac{1}{2(n-1)}\cdot\frac{\partial G(u,v)}{\partial v}
 =-\left[\lambda(n+1)-2\right]u+\left[\lambda(n-1)-2\right]v=0.
\end{array} \right .
\ee
The equations have only solution $(u,v)=(0,0)$, and
$G(0,0)=0$. We now divide the proof into four
cases.

\smallskip

\noindent {\bf Case 1}: $u=1$.\smallskip

Since $n\geq 3$ and $\frac{6n-2}{n^2+n}\leq\lambda$, by \eqref{eq7a}, it is clear that
$\lambda(n-1)-2>0$ and thus, the function $\psi (v)=G(1,v)$ on $[-1,1]$, where
\be\label{eq-ex2}
\psi (v)=[\lambda(n+1)^2-10n+2]-2(n-1)[\lambda(n+1)-2]v+(n-1)[\lambda(n-1)-2]v^2,
\ee
attains its maximum at $v=-1$ so that $\psi (v)\leq \psi (-1)$ for $v\in [-1, 1]$.

\smallskip

\noindent {\bf Case 2}: $u=-1$.

\smallskip

By similar reasoning as in Case 1, we can easily show that $G(-1,v)\leq G(-1,1)$ for $v\in [-1, 1]$.

\smallskip

\noindent {\bf Case 3}: $v=1$.

\smallskip

In this case, we regard $\phi (u)= G(u,1)$ as a function of $u$ on the interval $[-1,1]$. We have
\be\label{eq-ex3}
\phi (u)=[\lambda(n+1)^2-10n+2]u^2-2(n-1)[\lambda(n+1)-2]u+(n-1)[\lambda(n-1)-2]
\ee
so that
$$  \phi' (u)=2[\lambda(n+1)^2-10n+2]u-2(n-1)[\lambda(n+1)-2] ~\mbox{ and }~\phi'' (u) =2[\lambda(n+1)^2-10n+2].
$$
Solving the equation $ \phi' (u)=0$ gives  only solution
\be\label{eq-ex4}
u_0=\frac{(n-1)\left[\lambda(n+1)-2\right]}{\lambda(n+1)^2-10n+2}.
\ee
By assumption $u_0\leq -1$, where the inequality holds if and only if $\lambda=\frac{6n-2}{n^2+n}$.

If $\lambda=\frac{6n-2}{n^2+n}$, then $u_0=-1$ is the point of maximum of $\phi (u)$ on $[-1,1]$, since $\phi'' (u) <0$.

If $\lambda>\frac{6n-2}{n^2+n}$, then $u_0<-1$ and thus,
$$\phi (u)\leq \max\left\{ \phi(1), \phi(-1)\right\}=\phi (-1) ~\mbox{ for $u\in [-1,1]$}.
$$

\smallskip

\noindent {\bf Case 4}: $v=-1$.

\smallskip

By similar reasoning as in Case 3, we can easily see that
$$G(u,-1)\leq G(1,-1)~\mbox{ for $u\in [-1,1]$}.
$$
Finally, $G(-1,1)-G(0,0)=4\lambda n^2-16n+8>0,$ which is obvious, since $\lambda\geq\frac{6n-2}{n^2+n}$. The proof is completed.
\epf

%
%
%

\begin{lem}\label{G boundary bound 2}
Suppose that  $n\geq 3$ and  $0<\lambda\leq\frac{2n}{n^2-n+1}$. Then $G(u,v)\leq 0$ for $(u,v)\in R$, where
$G(u,v)$ is defined by \eqref{eq5}.
\end{lem}
\bpf
Clearly, we have the chain of inequalities
\be\label{eq-ex5}
\frac{2}{n+1}<\frac{2}{n}<\frac{2n}{n^2-n+1}<\frac{2}{n-1}<\frac{4n-2}{n^2}<\frac{5n-1}{n^2+n}<\frac{6n-2}{n^2+n}<\frac{8n}{(n+1)^2}
<\frac{10n-2}{(n+1)^2}.
\ee
\noindent {\bf Case 1}: $\lambda=\frac{2}{n+1}$.\smallskip

Substituting this value of $\lambda$ in \eqref{eq5}, we find that
$$G(u,v)=-4(2n-1)u^2-\frac{4(n-1)}{n+1}v^2\leq 0=G(0,0).
$$

\noindent {\bf Case 2}:  $\lambda=\frac{2n}{n^2-n+1}$.\smallskip

Again, substituting this value of $\lambda$ in \eqref{eq5}, we see that
$$G(u,v)=-\frac{2(n-1)}{n^2-n+1}\left[(2n-1)u+v\right]^2\leq 0=G(0,0).
$$

\noindent {\bf Case 3}:  $\lambda\notin \{\frac{2}{n+1},\, \frac{2n}{n^2-n+1}\}$.\smallskip

From the partial derivatives of $G(u,v)$ determined from \eqref{eq-ex1}, we obtain that the equations in \eqref{eq-ex1}
have only solution $(u,v)=(0,0)$, and $G(0,0)=0$. As before,  we need to divide the proof into four
subcases.

\medskip

\noindent {\bf Subcase 3(a)}: $u=1$. \smallskip

For this subcase, we consider the function $\psi (v)=G(1,v)$ defined by \eqref{eq-ex2}
and obtain that
$$\psi '(v)=- 2(n-1)\left[\lambda(n+1)-2\right]+2(n-1)\left[\lambda(n-1)-2\right]v.
$$
We see that $\psi '(v)=0$ yields the solution $v_0$, where
$$v_0=\frac{\lambda(n+1)-2}{\lambda(n-1)-2},  ~\psi (v_0)=\frac{8\left[2n-\lambda(n^2-n+1)\right]}{\lambda(n-1)-2}
~\mbox{and }~ \psi (-1)=4\lambda n^2-16n+8.
$$
A computation implies that $|v_0|\leq 1$ if and only if $0<\lambda\leq \frac{2}{n}$. Therefore, we have
$$G(1,v)\leq  \left \{ \begin{array}{rl}
 \ds \frac{8\left[2n-\lambda(n^2-n+1)\right]}{\lambda(n-1)-2} & \mbox{for $0<\lambda\leq \frac{2}{n}$}\\
\ds 4\lambda n^2-16n+8 & \mbox{for $\frac{2n}{n^2-n+1}>\lambda>\frac{2}{n} $.}
\end{array} \right .
$$

\noindent {\bf Subcase  3(b)}: $u=-1$.\smallskip

By similar reasoning as in Subcase~3(a), we can  easily see that the last inequality continues to hold
with $G(-1,v)$ instead of $G(1,v)$.

\noindent {\bf Subcase  3(c)}: $v=1$.\smallskip

For this subcase, we recall the function $\phi (u)= G(u,1)$ defined by \eqref{eq-ex3}. We observe that
$u_0$ defined by \eqref{eq-ex4} has the property
(by assumption) that $|u_0|<1$ for $0<\lambda<\frac{6n-2}{n^2+n}$ and so, $u_0$ is the point of maximum for
$\phi (u)$ on $[-1,1]$, which yields that
\be\label{eq-ex7}
\phi(u)\leq \phi (u_0)= \frac{8(n-1)\left[\lambda(n^2-n+1)-2n\right]}{10n-2-\lambda(n+1)^2}.
\ee

\noindent {\bf Subcase  3(d)}: $v=-1$.\smallskip

By similar reasoning as in Subcase 3(c), we can easily prove that the last inequality \eqref{eq-ex7} continues to hold in this case too.
Since $0<\lambda<\frac{2n}{n^2-n+1}$ and $n\geq 3$, we deduce that
$$\frac{8\left[2n-\lambda(n^2-n+1)\right]}{\lambda(n-1)-2}\leq0,~ 4\lambda n^2-16n+8<0,~\mbox{and }~
\frac{8(n-1)\left[\lambda(n^2-n+1)-2n\right]}{10n-2-\lambda(n+1)^2}\leq 0.
$$

Finally, the above facts imply the desired conclusion of the lemma.
\epf

\begin{lem}\label{F boundary bound}
Suppose that  $n\geq 3$,  $\frac{2n}{n^2-n+1}<\lambda<\frac{6n-2}{n^2+n}$ and $F(u,v) :=F_{n, \lambda}(u,v)$ is given by Lemma \ref{bound}. Then
$$\max\left\{F(u,v):\ (u,v)\in \partial R\right\}\leq A_{n,\lambda},
$$
where
\be\label{eq-ex6}
A_{n,\lambda} = \left \{
\begin{array}{rl}  \ds \frac{4(n-1)^2\left[\lambda(n-1)+1\right]}{8n-(n+1)^2\lambda}& \mbox{ for $\frac{2n}{n^2-n+1}<\lambda\leq \frac{5n-1}{n^2+n}$}\\[3mm]
 \ds  4\lambda n^2-12n+4& \mbox{ for $\frac{5n-1}{n^2+n}<\lambda<\frac{6n-2}{n^2+n}$}.
\end{array} \right .
\ee
\end{lem}
\bpf
We continue to use the chain of inequalities given by \eqref{eq-ex5} and as before, the proof is divided into four cases.\smallskip

\noindent {\bf Case 1}: $u=1$. \smallskip

For this case, $\lambda>\frac{2n}{n^2-n+1}>\frac{2}{n+1}$, the function $\Psi (v):=F(1,v)$ defined by
$$\Psi (v)=[(n+1)^2\lambda-8n]-2(n-1)[(n+1)\lambda-2]v+(n-1)^2\lambda v^2
$$
attains its maximum at $v=-1$. This observation yields that
$$\Psi (v)\leq \Psi (-1) =4\lambda n^2-12n+4.
$$

\noindent {\bf Case 2}: $u=-1$.\smallskip

By similar reasoning as in Case 1, we conclude that $F(-1,v)\leq  4\lambda n^2-12n+4.$

\noindent {\bf Case 3}: $v=1$.\smallskip

In this case, we need to consider the function $\Phi (u):=F(u,1)$ defined by
$$\Phi (u) =[(n+1)^2\lambda-8n]u^2-2(n-1)[(n+1)\lambda-2]u+(n-1)^2\lambda .
$$
Clearly,
the only solution $u_0$ to the equation $\Phi' (u)=0$ is given by
$$u_0=\frac{(n-1)\left[\lambda(n+1)-2\right]}{\lambda(n+1)^2-8n}.
$$
Moreover, $|u_0|\leq 1$ if and only if
$$[(n^2+n)\lambda-(5n-1)] \,[(n+1)\lambda-(3n+1)]\geq 0.
$$
Also, we have
$$\Phi (u_0)=F(u_0,1)= \frac{4(n-1)^2\left[\lambda(n-1)+1\right]}{8n-(n+1)^2\lambda} 
$$
and
$$\Phi (-1)=F(-1,1)=F(1,-1)=4\lambda n^2-12n+4.
$$
Since $\lambda<\frac{6n-2}{n^2+n}<\frac{3n+1}{n+1}$ and $\frac{5n-1}{n^2+n}<\frac{6n-2}{n^2+n}$, the above facts show that
$$|u_0| ~   \left \{
\begin{array}{rl}  \ds \leq 1 & \mbox{ for $\frac{2n}{n^2-n+1}<\lambda\leq \frac{5n-1}{n^2+n}$}\\[3mm]
 \ds  >1 &\mbox{ for $\frac{5n-1}{n^2+n}<\lambda<\frac{6n-2}{n^2+n}$}.
\end{array}
 \right .
$$
It follows that $F(u,1)\leq A_{n,\lambda},$
where $A_{n,\lambda}$ is given by \eqref{eq-ex6}.

 \noindent {\bf Case 4}. $v=-1$.\smallskip

Again, by similar reasoning as in Case 3, we can prove that $F(u,-1)\leq A_{n,\lambda}.$
Moreover, by a computation, for $\lambda<\frac{6n-2}{n^2+n}$, we have
$$\Phi (-1)-\Phi (u_0)=4\lambda n^2-12n+4-\frac{4(n-1)^2\left[\lambda(n-1)+1\right]}{8n-(n+1)^2\lambda}=-\frac{4\left[n(n+1)\lambda-(5n-1)\right]^2}{8n-(n+1)^2\lambda}\leq 0.
$$
The desired conclusion follows if we use the above facts and combine the four cases.
\epf

\begin{lem}\label{solution of equation}
Suppose that $n\geq 3$, $\frac{2n}{n^2-n+1}<\lambda<\frac{6n-2}{n^2+n}$, and $F(u,v) :=F_{n, \lambda}(u,v)$ is given by Lemma \ref{bound}.
Then the critical points of $F(u,v)$ are $(0, 0)$ and $(u, v)$, where $(u, v)$ satisfies
\be\label{eq12}
v^2=\frac{[\sqrt{n\lambda}+\sqrt{2-\lambda}]^2 [(n-1)\lambda-\sqrt{\lambda
n(2-\lambda)}]}{(n-1)^2\lambda^2},
\ee
\be\label{eq13}
uv=\frac{[(n+1)\lambda-2][\lambda(n-1)-\sqrt{\lambda
n(2-\lambda)}]}{(n-1)\lambda [\lambda(n-1)-2\sqrt{\lambda n(2-\lambda)}]},
\ee
and
\be\label{eq14}
u^2=\frac{[\lambda(n-1)-2\sqrt{\lambda n(2-\lambda)}+2][(n-1)\lambda-\sqrt{\lambda
n(2-\lambda)}]}{[\lambda(n-1)-2\sqrt{\lambda n(2-\lambda)}]^2}.
\ee
\end{lem}
\bpf 
The critical points of $F(u,v)$ are the solutions of
the equations
\be\label{eq7}
\frac{\partial F(u,v)}{\partial u}=0 ~\mbox{ and }~   \frac{\partial F(u,v)}{\partial v}=0,
\ee
which is equivalent to solving the pair of equations
\be\label{eq8} \left\{
\begin{array}{rl}  \ds \left[(n+1)^2\lambda-8n\right]u-(n-1)\left[(n+1)\lambda-2\right]v=2(n-1)u\frac{\sqrt{1-v^2}}{\sqrt{1-u^2}},\\
 \ds  -\left[(n+1)\lambda-2\right]u+(n-1)\lambda  v=2v\frac{\sqrt{1-u^2}}{\sqrt{1-v^2}}.
\end{array} \right.
\ee
It is obvious that $(u, v)=(0, 0)$ is a solution of \eqref{eq8}. We
now assume that $(u, v)\neq (0, 0)$. It follows from \eqref{eq8} that
\be\label{eq9}
a(n)\frac{u^2}{v^2}+b(n)\frac{u}{v}+c(n)=0,
\ee
where
\beqq
a(n)& =&  [(n+1)^2\lambda-8n ] [(n+1)\lambda-2 ]\\
b(n) &= &-2\lambda (n-1) [(n+1)^2\lambda-2(3n+1) ] ~\mbox{ and }\\
c(n)&=&\lambda(n-1)^2 [(n+1)\lambda-2].
\eeqq
We consider the discriminant $\triangle$ of the quadratic equation \eqref{eq9} in the variable $u/v$. By a computation, we have
$$\triangle= b^2(n)-4a(n)c(n)=64(n-1)^2\lambda n(2-\lambda)>0,
$$
since $n\geq 3$ and $\frac{2n}{n^2-n+1}<\lambda<\frac{6n-2}{n^2+n}$.
Thus, \eqref{eq9} has two solutions
$$\frac{u}{v}= (n-1)\left [\frac{\lambda  [(n+1)^2\lambda-2(3n+1)]\pm 4\sqrt{\lambda n(2-\lambda)}}
{[(n+1)^2\lambda-8n]\,[(n+1)\lambda-2]} \right ].
$$
In the case
$$u= (n-1)\left [\frac{\lambda  [(n+1)^2\lambda-2(3n+1)]-4\sqrt{\lambda n(2-\lambda)}}
{[(n+1)^2\lambda-8n]\,[(n+1)\lambda-2]} \right ] v,
$$
it follows from the second identity of \eqref{eq8} that
\beqq
2\sqrt{\frac{1-u^2}{1-v^2}}&=&(n-1) \left [\lambda-\frac{\lambda
[(n+1)^2\lambda-2(3n+1)]- 4\sqrt{\lambda
n(2-\lambda)}}{(n+1)^2\lambda-8n} \right ]\\
&=&-\frac{2(n-1)[\lambda(n-1)-2\sqrt{\lambda
n(2-\lambda)}]}{(n+1)^2\lambda-8n}.
\eeqq
The right of the last equation is less than zero, since
$$\lambda\left[(n+1)^2\lambda-8n\right]=\left[\lambda(n-1)+2\sqrt{\lambda
n(2-\lambda)}\right] \left[\lambda(n-1)-2\sqrt{\lambda n(2-\lambda)}\right]<0.
$$
This is clearly a contradiction.
In view of this observation,  we only need to consider the case
$$u= (n-1)\left [\frac{\lambda  [(n+1)^2\lambda-2(3n+1)]+4\sqrt{\lambda n(2-\lambda)}}
{[(n+1)^2\lambda-8n]\,[(n+1)\lambda-2]} \right ] v,
$$
which may be rewritten as
\be\label{eq10}
u=(n-1)\left [\frac{ \lambda [\lambda(n-1)-2\sqrt{\lambda n(2-\lambda)}+2]}{[(n+1)\lambda-2]
[\lambda(n-1)-2\sqrt{\lambda n(2-\lambda)}]}\right ]v.
\ee
It follows quickly from \eqref{eq10} and the second identity of \eqref{eq8} that
\be\label{eq11}
\sqrt{\frac{1-u^2}{1-v^2}}=\frac{(n-1)\lambda}{2\sqrt{\lambda n(2-\lambda)}-\lambda(n-1)}.
\ee
Substituting \eqref{eq10} in \eqref{eq11} yields
$$1- \frac{(n-1)^2\lambda^2[\lambda(n-1)-2\sqrt{\lambda
n(2-\lambda)}+2]^2}{[(n+1)\lambda-2]^2
[\lambda(n-1)-2\sqrt{\lambda
n(2-\lambda)}]^2}v^2=\frac{(n-1)^2\lambda^2}{[2\sqrt{\lambda
n(2-\lambda)}-\lambda(n-1)]^2}(1-v^2)
$$
which after simplification gives
$$v^2=\frac{[(n+1)\lambda-2]^2}{(n-1)^2\lambda^2}\cdot
\frac{(n-1)\lambda-\sqrt{\lambda n(2-\lambda)}}{[\sqrt{n\lambda}-\sqrt{2-\lambda}]^2}.
$$
This is the same as \eqref{eq12} and observe that the right side expression in the above form of $v^2$ is clearly positive, because
$$(n-1)^2 \lambda^2-\lambda n(2-\lambda)=\lambda[(n^2-n+1)\lambda-2n]>0.
$$
Using \eqref{eq12}, \eqref{eq10} easily gives \eqref{eq13} and  \eqref{eq14}. The proof of the lemma is complete.
\epf

\begin{lem}\label{solution of equation2}
Suppose that $F(u,v) :=F_{n, \lambda}(u,v)$ is given by Lemma \ref{bound}. Then we have the following:
\bee
\item [{\rm (1)}] If $n\geq 3$ and $\frac{3n+\sqrt{5n^2-4n}}{n^2+n}\leq\lambda<\frac{6n-2}{n^2+n},$
then the only critical point of $F(u,v)$ for $(u,v)\in (-1, 1)\times (-1, 1)$ is $(0, 0)$, and that
$F(0,0)=4(n-1).$

\item [{\rm (2)}] If $n\geq 3$ and $\frac{2n}{n^2-n+1}<\lambda<\frac{3n+\sqrt{5n^2-4n}}{n^2+n},$ then
there are three critical points for $F(u,v)$ on  $(-1, 1)\times (-1, 1)$, namely,
$(0, 0)$ and $(u_i, v_i)~ (i=1, 2)$ which satisfy \eqref{eq12}, \eqref{eq13} and \eqref{eq14}. Moreover,
\be\label{eq-lem8}
F(u_1,v_1)=F(u_2,v_2)=\frac{4\lambda(n-1)\left[\lambda(n^2+1)-4n\right]+16n(2-\lambda)\sqrt{n\lambda(2-\lambda)}}
{\lambda\left[8n-\lambda(n+1)^2\right]},
\ee
and $F(u_1,v_1)>F(0,0).$
\eee
\end{lem}
\bpf For $n\geq 3$ and $\frac{2n}{n^2-n+1}<\lambda<\frac{6n-2}{n^2+n}$, Lemma \ref{solution
of equation} and \eqref{eq11} yield $v^2>u^2.$ Equation \eqref{eq12} shows that $v^2<1$ if and only if
\be\label{eq15}
2\lambda\left[n\lambda-(n+1)\right]+\left[(n-1)\lambda-2\right]\sqrt{\lambda n(2-\lambda)}<0.
\ee
By assumption, we have $n\lambda-(n+1)<0$ which implies that \eqref{eq15} holds if $\frac{2n}{n^2-n+1}<\lambda\leq\frac{2}{n-1}$.

Thus, for $\frac{2}{n-1}<\lambda<\frac{6n-2}{n^2+n}$, \eqref{eq15} holds if and only if
$$A=:4\lambda\left[n\lambda-(n+1)\right]^2-\left[(n-1)\lambda-2\right]^2 n(2-\lambda)>0.
$$
By simplification, we find that
\beqq
A&=& n(n+1)^2\lambda^3-2n(n+1)(n+3)\lambda^2+4(3n^2+n+1)\lambda-8n\\
&=&\left[n(n+1)\lambda^2-6n\lambda+4\right]\left[(n+1)\lambda-2n\right]\\
&=&n(n+1)\left(\lambda-\frac{3n-\sqrt{5n^2-4n}}{n^2+n}\right)\left(\lambda-\frac{3n+\sqrt{5n^2-4n}}{n^2+n}\right)\left[(n+1)\lambda-2n\right].
\eeqq
Consequently, since  $\frac{2}{n-1}<\lambda<\frac{6n-2}{n^2+n}$ and
$$\frac{3n-\sqrt{5n^2-4n}}{n^2+n}<\frac{2}{n-1}<\frac{5n-1}{n^2+n}<\frac{3n+\sqrt{5n^2-4n}}{n^2+n}<\frac{6n-2}{n^2+n}<\frac{2n}{n+1},
$$
the above facts imply that $A>0$ if and only if
$$\frac{2}{n-1}<\lambda<\frac{3n+\sqrt{5n^2-4n}}{n^2+n}.
$$
We thus have shown that $1>v^2>u^2>0$ if
$$\frac{2n}{n^2-n+1}<\lambda<\frac{3n+\sqrt{5n^2-4n}}{n^2+n},
$$
whereas $v^2\geq 1$ if
$$\frac{3n+\sqrt{5n^2-4n}}{n^2+n}\leq\lambda<\frac{6n-2}{n^2+n}.
$$
Therefore, if
$$\frac{3n+\sqrt{5n^2-4n}}{n^2+n}\leq\lambda<\frac{6n-2}{n^2+n},
$$
there exists exactly one solution of the equation \eqref{eq7} for
$(u,v)\in (-1, 1)\times (-1, 1)$ and this solution is $(0, 0)$; if
$$\frac{2n}{n^2-n+1}<\lambda<\frac{3n+\sqrt{5n^2-4n}}{n^2+n},
$$
there exist three solutions of the equation \eqref{eq7}
for $(u,v)\in (-1, 1)\times (-1, 1)$, where one of them is $(0, 0)$ while the other two are $(u_i,
v_i) \,(i=1, 2)$ which satisfy \eqref{eq12}, \eqref{eq13} and \eqref{eq14}.

By a routine computation, the rest of the proof follows and we complete the proof.
 \epf

\begin{lem}\label{maximum of special cases}
Let $n\geq 3$ and $R=[-1,1]\times [-1,1]$. Then we have the following:
\bee
\item [{\rm (1)}] If $ \frac{3n+\sqrt{5n^2-4n}}{n^2+n}\leq\lambda<\frac{6n-2}{n^2+n},$ then
$F(u,v)\leq 4\lambda n^2-12n+4$ for $(u,v)\in R.$

\item [{\rm (2)}] If $ \frac{2n}{n^2-n+1}<\lambda<\frac{3n+\sqrt{5n^2-4n}}{n^2+n}$, then
$F(u,v)\leq F(u_1,v_1)$.
\eee
Here $F(u,v) :=F_{n, \lambda}(u,v)$ is given by Lemma \ref{bound} and $F(u_1,v_1)$ is given by \eqref{eq-lem8}.
\end{lem}
\bpf If $ \frac{3n+\sqrt{5n^2-4n}}{n^2+n}\leq\lambda<\frac{6n-2}{n^2+n},$
Lemmas \ref{F boundary bound} and \ref{solution of equation2} show
that for $(u,v)\in R$, we have
$$F(u,v)\leq \max\left\{4\lambda n^2-12n+4,\, 4(n-1)\right\}=4\lambda n^2-12n+4.
$$

If $ \frac{2n}{n^2-n+1}<\lambda<\frac{3n+\sqrt{5n^2-4n}}{n^2+n}$,
then we divide the proof into the following two cases.\medskip

\noindent {\bf Case 1}: $\frac{5n-1}{n^2+n}<\lambda<\frac{3n+\sqrt{5n^2-4n}}{n^2+n}$.
\smallskip

It follows from Lemmas \ref{F boundary bound} and \ref{solution of equation2}  that
$$F(u,v)\leq \max\left\{4\lambda n^2-12n+4,\,  F(u_1,v_1)\right\}=F(u_1,v_1),
$$
for

\vspace{8pt} $\ds \lambda\frac{8n-(n+1)^2\lambda}{4n}
\left[F(u_1,v_1)-(4\lambda n^2-12n+4)\right]$
\beqq
&=&\lambda\left[n(n+1)^2\lambda^2-2n(5n+3)\lambda+4(5n-1)\right]+4(2-\lambda)\sqrt{n\lambda(2-\lambda)}\\
&>& \frac{\sqrt{n\lambda(2-\lambda)}}{n-1}\left[n(n+1)^2\lambda^2-2(5n^2+5n-2)\lambda+(28n-12)\right]\\
&\geq&\frac{\sqrt{n\lambda(2-\lambda)}}{(n-1)n(n+1)^2}(3n^4-6n^3-n^2+8n-4)>0
 \eeqq
since $n\geq 3$,
$\frac{5n-1}{n^2+n}<\lambda<\frac{3n+\sqrt{5n^2-4n}}{n^2+n}$ and
$\lambda(n-1)>\sqrt{n\lambda(2-\lambda)}$.
\smallskip

\noindent {\bf Case 2}:
$\frac{2n}{n^2-n+1}<\lambda\leq\frac{5n-1}{n^2+n}$.
\smallskip

Lemmas \ref{F boundary bound} and \ref{solution of equation2} imply
that
$$F(u,v)\leq \max\left\{\frac{4(n-1)^2\left[\lambda(n-1)+1\right]}{8n-(n+1)^2\lambda}, \,  F(u_1,v_1)\right\} =F(u_1,v_1),
$$
for

\vspace{8pt}

$\ds \lambda\frac{8n-(n+1)^2\lambda}{4}
\left[F(u_1,v_1)-\frac{4(n-1)^2\left[\lambda(n-1)+1\right]}{8n-(n+1)^2\lambda}\right]$
\beqq
&=&\lambda(n-1)(2n\lambda-5n+1)+4n(2-\lambda)\sqrt{n\lambda(2-\lambda)}\\
&>& \sqrt{n\lambda(2-\lambda)}(2n\lambda-5n+1+8n-4n\lambda)\\
&=&2n\sqrt{n\lambda(2-\lambda)}\left(\frac{3n+1}{2n}-\lambda\right)>0,
\eeqq
since $n\geq 3$, $\frac{2n}{n^2-n+1}<\lambda\leq\frac{5n-1}{n^2+n}$ and
$\lambda(n-1)>\sqrt{n\lambda(2-\lambda)}$. The proof is completed.
\epf

\section{Proof of Theorem \ref{main theorem}}\label{sec4}

Suppose that $f\in {\mathcal C}$ and $f(z)=z+\sum_{n=2}^{\infty}a_{n}z^{n}$.
Since $|\lambda a^2_n-a_{2n-1}|$ is invariant under rotations, we consider instead the problem of maximizing the functional
${\rm Re }\,( \lambda a^2_n-a_{2n-1})$. Moreover, we find that
$${\rm Re }\,( \lambda a^2_n-a_{2n-1})=\lambda \left({\rm Re\, }
a_{n}\right)^2-\lambda \left({\rm Im\, } a_{n}\right)^2-{\rm Re\, } a_{2n-1}\leq J(f)
=\lambda \left({\rm Re\, } a_{n}\right)^2 -{\rm Re\, } a_{2n-1}.
$$
It follows from Lemmas \ref{bound}, \ref{real function bound}, \ref{G boundary bound 1} and \ref{maximum of special cases}(1) that for $n\geq 3$,
$$4J(f)\leq 4n + F(u,v) \leq 4n + (4\lambda n^2-12n+4)=4[\lambda n^2-(2n-1)]
$$
and the desired conclusion of Theorem \ref{main theorem}(1)  follows.

With the same reasoning as above, Lemmas \ref{bound} and \ref{maximum of special cases}(2) give Theorem \ref{main theorem}(2).

Finally, by the inequality \eqref{eq4}, Lemmas \ref{bound} and \ref{G boundary bound 2}, we obtain that, for $0<\lambda\leq\frac{2n}{n^2-n+1}$,
$$4J(f)\leq 4n +4(n-1)+G(u,v)\leq 4(2n-1)
$$
and Theorem \ref{main theorem}(3) follows.
\hfill $\Box$

\subsection*{Acknowledgements}
This work was completed during the visit of the first author to
Syracuse University. She thanks the university for its hospitality.
The visit and the research of the first author was supported by CSC
of China (No. 201308430274). The research was also supported by NSF
of China (No. 11201130 and No. 11571216), Hunan Provincial Natural
Science Foundation of China (No. 14JJ1012) and construct program of
the key discipline in Hunan province. The second author is on leave from
IIT Madras.

\end{document}